\documentclass[a4paper]{article}

\usepackage{graphics}
\usepackage{amsmath}
\usepackage{amssymb}
\usepackage{amscd}

\newcommand{\tr}[1][]{\ensuremath{\stackrel{#1}{\longrightarrow}}}
\newcommand{\proof}{\textbf{Proof : }}
\newcommand{\cqfd}{\hfill \fbox{} \vskip 0.2cm}
\newcommand{\edge}[1]{\ensuremath{\lbrace #1 \rbrace}}
\newcommand{\sh}[1]{\ensuremath{s(#1)}}

\newcommand{\fig}[3][1]{
\begin{figure}[!h]
\centerline{\scalebox{#1}{\includegraphics{#2.eps}}}
\caption{\label{fig_#2} \em #3} \end{figure}}

\newtheorem{theo}{Theorem}

\newtheorem{lem}{Lemma}

\begin{document}

\begin{center}
%{\huge Equivalence of Edge Firing Games and Distributive Lattices}\\
{\huge Coding Distributive Lattices\\ with Edge Firing Games}\\
\vskip 0.2cm
{\Large Matthieu Latapy and Cl\'emence Magnien\\
\vskip 0.2cm
\large \textsc{liafa} -- Universit\'e Paris 7\\
 2 place Jussieu, 75005 Paris.\\
(latapy,magnien)@liafa.jussieu.fr}
\vskip 0.2cm
\end{center}

\noindent
\textbf{Abstract:}
In this note, we show that any distributive lattice is isomorphic to
the set of reachable configurations of an Edge Firing Game.
Together with the result of James Propp, saying that the set of reachable
configurations of any Edge Firing
Game is always a distributive lattice, this shows that the two concepts
are equivalent.

\vskip 0.2cm

\noindent
\textbf{Keywords:}
Edge Firing Game, Source Reversal Game, Orientations of Graphs,
Distributive Lattice,
Discrete Dynamical Model, Chip Firing Game.

\vskip 0.5cm

\section{Background}

The Edge Firing Game (EFG) has been introduced in the context of graph
flow studies \cite{MKM78}, and has been re-introduced in various
occasions, as for example in order theory \cite{Pre86a,Pre86b}. Since then,
it has been widely studied, mainly from a combinatorial point of view
\cite{GP97,Pro93}.

Given an undirected graph $G=(V,E)$, an \emph{orientation} of $G$ is
a directed graph $G'=(V',E')$ such that $V'=V$, $(v,v') \in E'$
implies $\lbrace v, v' \rbrace \in E$, and $\lbrace v, v' \rbrace \in E$
implies either $(v,v') \in E'$ or $(v',v) \in E'$. An EFG is defined
by a connected undirected graph $G$ with a distinguished vertex $s$,
called the \emph{sink}, and an orientation $G_0$ of this graph. A
\emph{configuration} of the game is an orientation of $G$ and $G_0$ is
called the \emph{initial configuration}. The game is played with
respect to the following rule: one can transform a configuration $C$
into the configuration $C'$, which is denoted by $C \tr C'$, if there
is in $C$ a vertex $\nu \not= s$ which has
only incoming edges (and no outgoing edges),
and if $C'$ is obtained from $C$ by reversing all these edges, i.e.
by replacing each edge $(v,\nu)$ by $(\nu,v)$. 
We call this \emph{firing} $\nu$.
Notice that, if $C$
is an orientation of $G$ then $C'$ is another orientation of $G$.
If we iterate this rule starting from the initial configuration,
we obtain a set of reachable configurations, between which the evolution
rule $\tr$ defines a relation, called the \emph{successor relation}.
The set of reachable configurations,
together with the successor relation,
is called the \emph{configuration space} of
the game. All the orientations of $G$ which are in the configuration
space share some flow properties, which are detailed in \cite{MKM78}.

We recall that 
an \emph{ordered set} (or \emph{partially ordered set})
is a set equipped with a binary relation $\le$ which is
reflexive ($x \le x$), transitive ($x \le y$ and
$y \le z$ implies $x \le z$) and antisymetric ($x \le y$ and $y \le x$
implies $x=y$).
Given two elements $x$ and $y$ of an ordered set, we say that $x$
\emph{covers} $y$ (or $y$ is covered by $x$) 
and we write $x\succ y$ (or $y\prec x$)
if $x > y$ and
$x\ge z > y$ implies $z=x$.
A \emph{maximal} element in an ordered set is an element
which is smaller than no element.
A {\em lattice} is an ordered set such that
any two elements $a$ and $b$ have a least
upper bound (called the {\em join} of $a$ and $b$ and denoted by $a \vee b$)
and a greatest lower bound (called the {\em meet} of $a$ and $b$ and denoted by
$a \wedge b$). The element $a \vee b$ is the smallest element among the
elements greater than both $a$ and $b$. The element $a \wedge b$ is
defined dually. A  lattice
is \emph{distributive} if for all $a$, $b$ and $c$:
$(a \vee   b)\wedge (a \vee   c) = a \vee   (b \wedge c)$ and 
$(a \wedge b)\vee   (a \wedge c) = a \wedge (b \vee   c)$.
A distributive lattice is a strongly structured set, and many
general results, for example efficient coding and algorithms,
are known about such sets. For more details,
see for example \cite{DP90}.

One can prove that there can be no loop in the configuration space
of any EFG, therefore the successor relation induces an order over
the configurations, which is the transitive and reflexive closure of
this relation: $C \ge C'$ if and only if there exists a sequence of
firings which transforms $C$ into $C'$.
James Propp \cite{Pro93} proved much more on this
relation: it gives the distributive lattice
structure to the configuration space. The aim of this note is to prove the converse: given a
distributive lattice $L$, one can construct an EFG such that
its configuration space is isomorphic to $L$.

\section{The result}

In this section, 
we prove that any distributive lattice is isomorphic to
the configuration space of an EFG. 
The proof uses the theorem of Birkhoff for representation of
distributive lattices \cite{Bir33}. 
We first state this theorem and the
definitions needed.
After this, we restrict ourselves to the EFGs which satisfy the following
property: during any sequence of possible firings, each vertex is
fired at most once. Such an EFG is said to be \emph{simple}. We prove
that any distributive lattice is isomorphic to the configuration space
of such an EFG. Since, from \cite{Pro93}, the configuration space of
any EFG is a distributive lattice, we obtain as a corollary
that any EFG is equivalent
to a simple one.

Given a lattice $L$, an element $a$ of $L$ is a \emph{join-irreducible} if there is a
unique element $b$ in $L$ such that 
$b\prec a$.
If one considers the set $J$ of all the join-irreducibles
of a lattice $L$, then
a natural order is induced by $L$ over $J$: $a \le b$
in $J$ if and only if $a \le b$
in $L$. See Figure~\ref{fig_ex_ordreJ} for an example.

\fig[0.4]{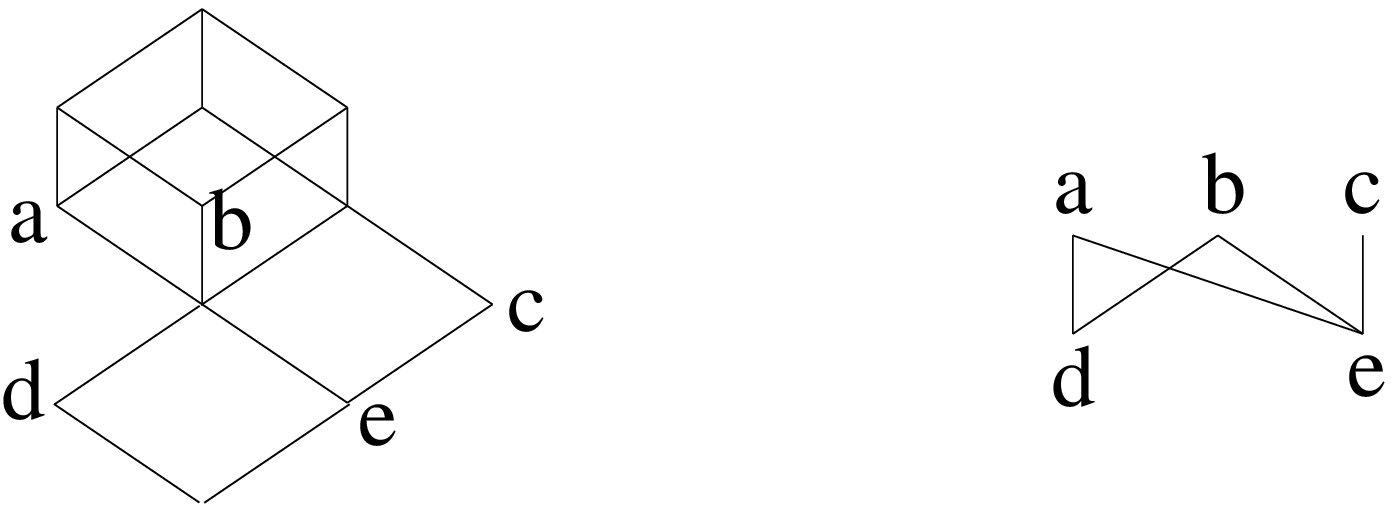}{A distributive lattice, and the order induced over
its join-irreducibles.}

Given an ordered set $O$, a \emph{filter} $F$ of $O$ is a subset of $O$
such that $a \in F$,
$b \in O$, and $a \le b$ imply $b \in F$. In other words, a
subset of $O$ is a filter of $O$ if and only if it is upper closed.
The set of all the filters of an ordered set can itself be ordered
by reverse inclusion:
$F \le F'$ if and only if $F' \subseteq F$. We can now state the famous
theorem from Birkhoff:

\begin{theo} \cite{Bir33}
\label{th_bir}
Any distributive lattice is isomorphic to the set of the filters of the
order induced over its join-irreducibles, ordered by reverse inclusion.
\end{theo}

In the sequel, we will only consider simple EFGs. Indeed, we will
see that any EFG is equivalent to such an EFG.
We begin by giving a lemma about simple EFGs, then we will prove our result.

\begin{lem}
\label{lem_shot_set}
Let $E$ be a simple EFG.
Let $C$ and $C'$ be two configurations of $E$ such that $C\ge C'$.
If $\sigma$ and $\sigma'$ are two sequences of firings
which transform $C$ into $C'$, then the set of fired vertices
during $\sigma$ is the same as during $\sigma'$.
\end{lem}
\proof
Let $S$ and $S'$ be the sets of vertices fired respectively during
the sequences $\sigma$ and $\sigma'$.
When we go from $C$ to $C'$ following $\sigma$, the
orientation of the edges of the graph vary as follows:
the edges between vertices not in $S$ are not affected,
the edges incident to a vertex in $S$ are reversed (they were directed
      towards $S$ in $C$, and are directed away from $S$ in $C'$),
and the edges between two vertices in $S$ have the same
      orientation in $C$ and $C'$
      (they are reversed twice).
From this, one can easily see that, if $S'$ is not equal to $S$, then the
sequence $\sigma'$ cannot lead to the same configuration as $\sigma$, which
proves the claim.
\cqfd

Notice that this lemma remains true if we consider EFGs which are not simple.
This can be shown with a proof similar to the one concerning Chip Firing
Games (see \cite{BLS91,LP00}). However, this proof would be
much more complicated than the one given here, which is sufficient to make
the note self-contained.

This lemma makes it possible to define the \emph{shot-set} of any
configuration $C$: \sh{C} is the set of vertices fired
to reach $C$ from the initial configuration.
We can now prove our result:

\begin{theo}
\label{th_efg_distr}
Any distributive lattice is isomorphic to the configuration space of
a simple Edge Firing Game.
\end{theo}
\proof
Let $L$ be a distributive lattice and $J$ be the set of its join-irreducibles
with the order induced by $L$.
Let $G=(V,E)$ be the undirected graph
defined as follows: $V = J \cup \lbrace \bot \rbrace$ where $\bot \not\in J$,
and $E = \lbrace \lbrace j, j' \rbrace\ |\ j \prec j' \mbox{ or }
 j' \prec j \mbox{ in } J \rbrace \ \cup\ \lbrace \lbrace \bot, j \rbrace\ |\ j
\mbox{ is a maximal element of } J \rbrace$.
We also define $G_0=(V,E_0)$
 in the following way:
$E_0\ =\ \lbrace (j,j')\ |\ 
j \prec j' \mbox{ in } J \rbrace
\ \cup\ \lbrace (\bot,j)\ |\ j
\mbox{ is a maximal element of } J \rbrace$. Notice that $G_0$ is an
orientation of $G$.

We claim that the configuration space of the EFG played on $G$,
with initial configuration $G_0$ and with sink $\bot$, is isomorphic to $L$.
To prove this, we will first show that the EFG is simple (\emph{i.e.}
each vertex is fired at most once), and then we will show that
the configuration space of the game and
the lattice of the filters of $J$ are isomorphic.

We first show by induction that each vertex is fired
at most once.
Let $j$ be a maximal element. We will
show that $j$ is fired only once.
It can be fired in $G_0$
since, by definition, all the edges incident to $j$ are directed towards $j$.
Once $j$ has been fired, the edge \edge{\bot,j} is directed
from $j$ to $\bot$, and since $\bot$ cannot be fired,
it will never be directed towards $j$ again.
Therefore $j$ can never be fired again.
Let now $j$ be an element covered by another element
$j'$ which can be fired only once.
We will show 
that $j$ can be fired at most once. It
can be fired only after $j'$
(the edge \edge{j,j'} must be directed towards $j$).
After the firing of $j$, the edge $\{j,j'\}$ is directed
towards $j'$, and since $j'$ cannot be fired again, the edge
will never be directed towards $j$ again, and $j$ cannot be fired again. 
Finally, the EFG is simple and thus, from Lemma~\ref{lem_shot_set},
we can consider the shot-sets of its configurations.

Now we prove that $s$, the function which associates its shot-set to a
configuration, is an isomorphism between the set of configurations of
the game and 
the set of filters of $J$.
To do this, we will first show that, for any configuration
$C$, \sh{C} is a filter, and then that, if $C_1,\ldots,C_k$ are the
configurations covered by $C$, then the
filters $\sh{C_1},\ldots,\sh{C_k}$ are exactly the ones covered by \sh{C}.
We will show these two steps simultaneously by induction.
To prove these two steps, we also need to prove that the vertices which
can be fired in configuration $C$ are exactly the maximal elements of
$J\setminus \sh{C}$. Therefore we will also prove this in the induction.
The base case is the initial configuration $G_0$:
 \sh{G_0} is the empty
set, which is the maximal filter for reverse inclusion.
The vertices which can be fired in $G_0$ are, by definition, the
maximal elements of $J$,
which are the maximal elements of $J\setminus\sh{G_0}$, as expected.
The induction step is as follows:
if $C$ is a configuration such that \sh{C} is a filter,
and such that the vertices that can be fired are the maximal elements of $J\setminus\sh{C}$,
we will show the following:
if $C_1,\ldots,C_k$ are the configurations covered by $C$, then the
sets $\sh{C_1},\ldots,\sh{C_k}$ are some filters, and they are exactly the
filters covered by \sh{C}.
We will prove simultaneously that the vertices which can be fired in configuration
$C_i$ are the maximal elements of $J\setminus\sh{C_i}$.
Let $C$ be any configuration such that \sh{C} is a filter,
and such that the vertices that can be fired are the maximal elements of $J\setminus\sh{C}$.
First notice that, in $J$,
the filters covered by a filter $F$ are exactly the sets $F\cup\{x\}$ for
all $x$ maximal in $J\setminus  F$.
Let $C_1,\ldots,C_k$ be the configurations covered by $C$.
Each $C_i$ is obtained from $C$ by the firing of a vertex $j_i$,
therefore $\sh{C_i}=\sh{C}\cup \{j_i\}$.
Since, by induction hypothesis, the $j_i$ are the maximal elements of $J\setminus\sh{C}$,
 the sets \sh{C_i} are exactly the filters
covered by \sh{C}, as expected.
To complete the induction step we show that, for all $i$, the vertices
that can be fired in configuration $C_i$ are exactly the maximal elements
of $J\setminus \sh{C_i}$.
Any maximal element $x$ of $J\setminus \sh{C_i}$ can be fired because
all the elements greater than $x$ (including the elements covering $x$) 
have been fired. Therefore all the edges \edge{x,x'} with $x\prec x'$
are now directed towards $x$ (the edges \edge{x,x'} with $x\succ x'$ have
not been turned and are still directed towards $x$).
Now we show that no vertex that is not maximal in $J\setminus \sh{C_i}$
can be fired:
 if a vertex $v$ is not maximal in
$J\setminus \sh{C_i}$, then at least one element $u$ that covers $v$ has not been
fired, and the corresponding edge is directed away from $v$.
This completes the induction step.

We have shown that, for any configuration $C$, \sh{C} is a filter,
and that the configurations covered by $C$ correspond exactly to the filters
covered by \sh{C}.
Therefore the configuration space of the EFG is isomorphic to the
lattice of the filters of the order induced over $J$, which by
Theorem~\ref{th_bir} is isomorphic to $L$.
\cqfd

\noindent
This theorem is illustrated in Figure~\ref{fig_ex_efg_distr}.
Notice that, since the configuration space of any EFG is a distributive
lattice, we obtain as a corollary that any EFG is equivalent to a simple
one. Therefore, we finally have that distributive lattices, EFGs
and simple EFGs are equivalent in terms of configuration spaces.
Notice also that the proof is constructive, which gives a way to
transform any EFG into a simple one.

\fig[0.5]{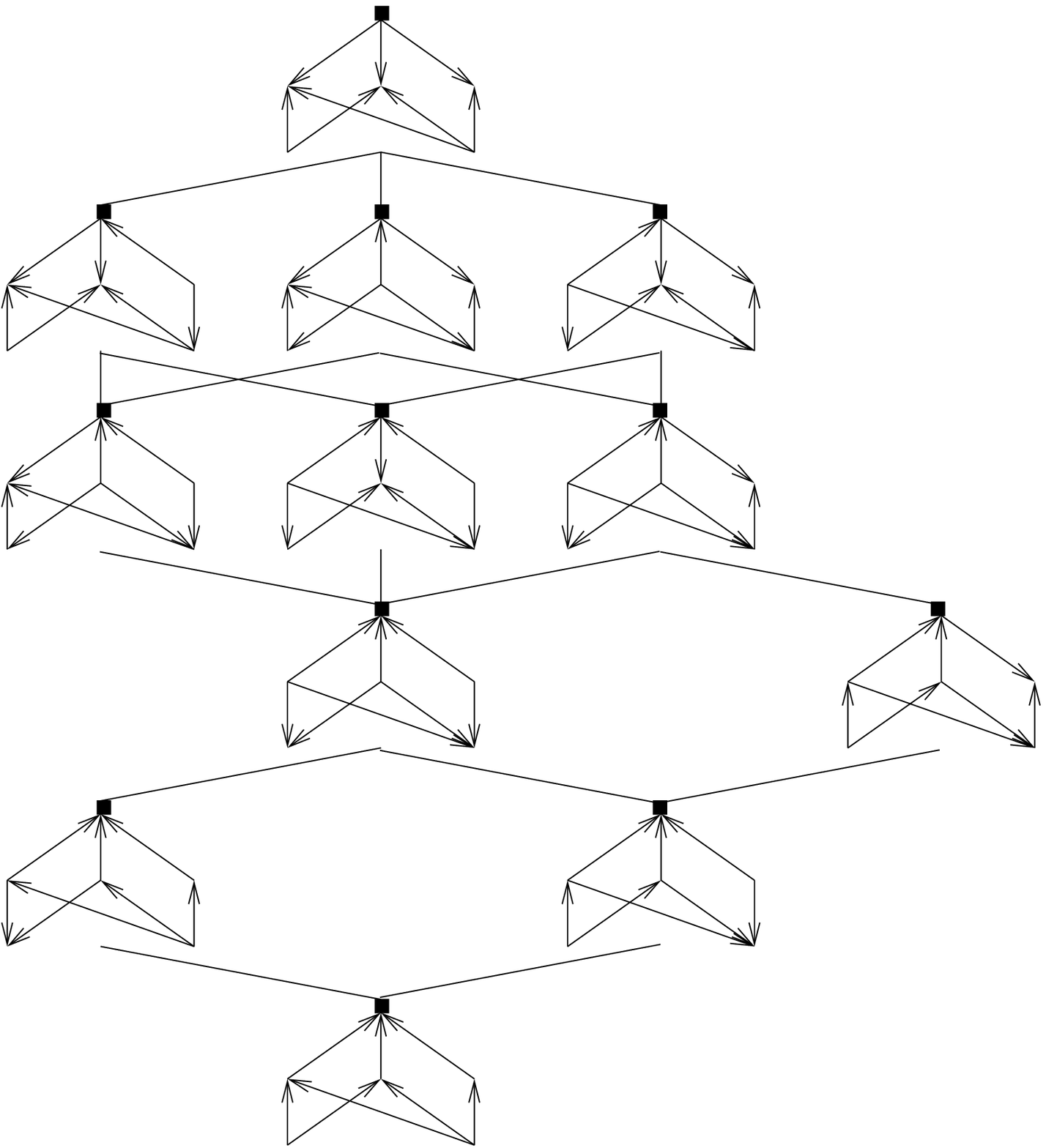}{The
EFG obtained by the proof of
Theorem~\ref{th_efg_distr}: its configurations space is isomorphic to
the distributive lattice shown in Figure~\ref{fig_ex_ordreJ}.   
The sink of the EFG is marked with a black square.}

\section{Discussion and perspectives}

The study on what kind of structures can be generated using
discrete dynamical models is an active area of research, and it is
now clear that lattices and distributive lattices often appear in this context.
In particular, two important models which appear in physics, combinatorics,
computer science and
social science, the Chip Firing Game (CFG) and the Abelian Sandpile Model (ASM),
have the following two properties: any distributive lattice can be obtained
as the configuration space of these models, and all the configuration
spaces one can obtain are Upper Locally Distributive (ULD) lattices.
These two classes of lattices are strongly structured and very close, therefore
we have a precise idea of what kind of configurations spaces are obtained
with these models. However, no exact characterization is known,
despite the fact that we know that they are all different.
The result we prove in this note, together with the result from
James Propp \cite{Pro93}, shows that there is an
equivalence between the distributive
lattices and EFGs. Therefore,
the EFG is strictly less powerful, in terms of the
obtainable configuration spaces, than the ASM and the CFG.

\small
\bibliographystyle{plain}
\bibliography{../Bib/bib.bib}

\end{document}